\documentclass[preprint,12pt]{article}   

\usepackage{graphicx}

\usepackage{amssymb}
\usepackage{amsmath}
\usepackage{amsthm}
\usepackage{bbm}
\usepackage{hyperref}
\usepackage{tabularx}
\usepackage{color}
\usepackage{cases}
\usepackage{stmaryrd}
\usepackage{geometry}
\newcommand{\R}{{\mathbb R}}
\newcommand{\psd}{p.s.d.}
\newcommand{\pd}{p.d.}
\newcommand{\psdsymb}{\succeq}
\newcommand{\pdsymb}{\succ}

\newcommand{\var}{\mathrm{Var}}
\newcommand{\cov}{\mathrm{Cov}}
\newcommand{\diag}{\mathrm{diag}}

\newcommand{\mean}[1]{\overline{#1}}
\newcommand{\avemap}{\phi}
\newcommand{\fillmap}{\psi}
\newcommand{\intInt}[1]{\llbracket 1, #1 \rrbracket}
\newcommand{\groupGeqTwo}[1]{\mathcal{G}^{\geq 2}_{#1}}
\newcommand{\mat}[1]{\mathbb{M}_{#1}}
\newcommand{\blockClass}[1]{\mathcal{B}_{#1}}
\newcommand{\blockClassG}[1]{\mathcal{C}_{#1}}
\newcommand{\blockClassI}[1]{\mathcal{B}^{\textrm{bg-iso}} _{#1}}
\newcommand{\compound}[2]{\Gamma^{\mathrm{CS}}_{#1} \left(#2\right)}
\newcommand{\symMat}[1]{\mathcal{S}_{#1}}

\newcommand{\one}[1]{\mathbbm{1}_{#1}}
\newcommand{\m}{}
\newcommand{\bs}{} %\{boldsymbol}

\usepackage{authblk}   

\title{On the validity of parametric block correlation matrices 
        with constant within and between group correlations}

\author[1]{O. Roustant}
\author[2]{Y. Deville}
\affil[1]{{\small Mines Saint-\'{E}tienne, UMR CNRS 6158, LIMOS, F--42023 Saint-\'{E}tienne, France}}
\affil[2]{{\small AlpeStat, Chamb\'ery, France}}
\begin{document}

\maketitle

\begin{abstract}
We consider the set $\blockClass{p}$ of parametric block correlation matrices with $p$ blocks of various (and possibly different) sizes,
whose diagonal blocks are compound symmetry (CS) correlation matrices and off-diagonal blocks are constant matrices. 
Such matrices appear in probabilistic models on categorical data, when the levels are partitioned in $p$ groups,
assuming a constant correlation within a group %(within group correlation) 
and a constant correlation for each pair of groups.\\ %(between group correlations).\\
We obtain two necessary and sufficient conditions for positive definiteness of elements of $\blockClass{p}$.
Firstly we consider the block average map $\phi$, consisting in replacing a block by its mean value.
We prove that for any $A \in \blockClass{p}$, $A \pdsymb 0$ if and only if $\avemap(A) \pdsymb 0$.
Hence it is equivalent to check the validity of the covariance matrix of group means,
which only depends on the number of groups and not on their sizes.
This theorem can be extended to a wider set of block matrices.
Secondly, we consider the subset $\blockClassI{p}$ of $\blockClass{p}$ 
for which the between group correlation is the same for all pairs of groups. 
Positive definiteness then comes down to find the positive definite interval of a matrix pencil on $\symMat{p}$.
We obtain a simple characterization by localizing the roots of the determinant with within group correlation values.
\end{abstract}
%\end{frontmatter}

%\keyword{titi, toto, tutu}

\section{Introduction}
We consider the set $\blockClass{p}$ of real parametric block correlation matrices
$$ A = \begin{pmatrix}
    A_{1,1}  & \dots & A_{1,p} \\
    \vdots & & \vdots \\
	A_{p,1}  & \dots & A_{p,p}  \\
    \end{pmatrix}$$
where diagonal blocs $A_{k,k}$ are compound symmetry matrices of size $n_k$ ($k \in \llbracket 1, p \rrbracket$)
and off-diagonal blocks $A_{k, \ell}$ are constant matrices of size $n_k \times n_\ell$:
$$ A_{k,k} = \underset{n_k}{\underbrace{\begin{pmatrix}
    1  & b_k & \dots &  b_k \\
    b_k & \ddots & \ddots & \vdots \\
    \vdots & \ddots & \ddots &  b_k\\
    b_k & \dots & b_k & 1\\
    \end{pmatrix}}}, 
    \qquad 
   A_{k,\ell} = A_{\ell, k}^\top = c_{k,l} \one{n_k} \one{n_\ell}^\top \quad (k < \ell), $$
All parameters $b_k, c_{k,\ell}$ are assumed to be in $]-1, 1[$.
When $n_k = 1$ we set by convention $b_k = 0$, 
and denote by $\groupGeqTwo{p}$ the set of $k \in \intInt{p}$ such that $n_k \geq 2$.\\

Such matrices appear in probabilistic models on categorical data, when the levels are partitioned 
in $p$ groups $G_1, \dots, G_p$ of sizes $n_1, \dots, n_p$,
assuming a constant correlation $b_k$ within the group $G_k$ %($k \in \llbracket 1, p \rrbracket$)
and a constant correlation $c_{k,l}$ for each pair of groups $(G_k, G_\ell)$.
Without loss of generality, we assume that $G_1$ corresponds to the first $n_1$ levels, 
$G_2$ to the $n_2$ next ones, and so on.
This results in a parameterization of the correlation matrix (see e.g. \cite{Qian_Wu_Wu_2008}) 
involving a small number of parameters.
In order to further reduce the number of parameters, one may often consider the subset 
$$\blockClassI{p} = \{ A \in \blockClass{p} \textrm{ s.t. for all } k,l: c_{k,l} = c \}$$ 
such that the between group correlation has a common value over pairs of groups.
Even more parsimonious models are obtained when a common value is  also assumed for within group correlations.\\

There are some common points with statistical models involving groups of variables, 
such as linear mixed effects models or hierarchical models.
When $p = 1$, $A$ is simply the CS correlation matrix used in random intercept model.
When the blocks have the same size, a common within group correlation and a common between group correlations, 
$A$ is a block CS matrix met in equicorrelated models. In that case, it is also connected to linear models with doubly exchangeable distributed errors \cite{Roy_Fonseca_2012}.\\ 

We aim at answering the question: For which values of $b_k, c_{k, \ell} \in ]-1, 1[$  ($1 \leq k < l \leq p$) 
is the matrix $A$ positive definite? 
In the simple case where $A$ has a block CS structure, as cited in the previous paragraph, 
the necessary and sufficient condition for positive definiteness is known 
(see e.g. \cite{Ritter_Gallegos_2002}, Lemma 4.3.). 
There is also a connection to matrix pencils, corresponding to matrices parameterized linearly by a single parameter. 
However, in general, the result seems not to be known.\\  

The paper is structured as follows. Section~\ref{sec:preliminaries} fixes notations and gives prerequisites.
The main results are given in Section~\ref{sec:main}, followed by corollaries in Section~\ref{sec:cor}.
An illustration is provided in Section~\ref{sec:num}. 
%Finally, Section~\ref{sec:extension} gives an extension of the first main result to a larger set of block matrices.
For the sake of readability, the proofs are gathered in a separate section (Section~\ref{sec:proofs}).   

\section{Preliminaries} \label{sec:preliminaries}
\subsection{Additional notations and basics}
\begin{itemize}
\item \emph{Positive definiteness.} If $A \in \mat{n}$, recall that $A$ is \emph{positive semidefinite (\psd)}, 
or simply \emph{positive}, if $\bs{x}^\top A \bs{x} \geq 0$ for all $\bs{x} \in \mathbb{R}^n$. We denote: $A \psdsymb 0$.
Similarly $A$ is \emph{positive definite (\pd)} if $\bs{x}^\top A \bs{x} > 0$ for all non-zero $\bs{x} \in \mathbb{R}^n$. 
We denote: $A \pdsymb 0$.
\item \emph{Matrix of ones}. %$\one{p}$ is the vector $p \times 1$ of ones, 
$J_p = \one{p} \one{p}^\top $ denotes the $p \times p$ matrix of ones. 
More generally, $\m{J}_{p,q} := \one{p} \one{q}^\top$.
\item \emph{Compound symmetry} (CS) matrices. 
$\compound{p}{v, c} = v I_p + c(J_p - I_p)$ is the CS covariance matrix. 
The CS correlation matrix is $\compound{p}{1, \rho}$, 
and for $v \neq 0$, $\compound{p}{v, c} = v \compound{p}{1, c/v}$ %\frac{c}{v}}$. 
It is well-known that 
\begin{equation} \label{eq:CompoudCarac} 
\compound{p}{v, c} \pdsymb 0 \quad \iff \quad - \frac{v}{p-1} < c < v
\end{equation}
For instance denoting $ b = c/v$, %\frac{c}{v}$, 
one can see that the eigenvalues of $\compound{p}{v, c}$ 
are  $1+(p-1)b$ with multiplicity $1$ (eigenvector $\one{p}$) 
and $1 - b$ with multiplicity $p-1$ (eigenspace: $\one{p}^\perp$).
\end{itemize}

\subsection{The block average map}
We denote by $\phi$ the \emph{block average map} on $\blockClass{p}$,
consisting in replacing a block by its mean value:
$$ A \in  \blockClass{p} \mapsto [\avemap(A)]_{k,l} = \frac{1}{n_k n_\ell} \sum_{i \in G_k, j \in G_\ell}{A_{i,j}}.$$
$\phi$ is an example of \emph{positive linear map} 
in the sense that if $A \psdsymb 0$ then $\avemap(A) \psdsymb 0$. 
Indeed, $\avemap(A) = V^\top A V$ where $V$ is the $p \times n$ matrix defined by
$$ \m{V} := 
\begin{pmatrix} 
  \frac{1}{n_1}\one{n_1} & 0 & \dots & 0 \\
  0 & \frac{1}{n_2}\one{n_2} & \ddots & \vdots \\
  \vdots & \ddots & \ddots & 0\\
  0 & \dots & 0 & \frac{1}{n_p}\one{n_p} \\
\end{pmatrix}. 
$$
and we retrieve one of the typical cases presented in \cite{Bhatia_2007_bookPD} (Example 2.2.1.(vi)).
This also can be viewed with a probabilistic interpretation. 
If $A \psdsymb 0$, it is a covariance matrix of some random vector $\bs{Y} = (Y_1, \dots, Y_n)$.
Then $\avemap(A)$ is the covariance matrix of the group means $(\overline{Y}_1, \dots, \overline{Y}_p)$ 
where $\overline{Y}_k = \frac{1}{n_k} \sum_{i \in G_k} Y_{i}$.\\
Finally, we give the explicit form of $\avemap(A)$ when $A \in \blockClass{p}$.  
Define for $k \in \llbracket 1, p \rrbracket$:
\begin{equation} \label{eq:alpha}
\alpha_k = \frac{1}{n_k} + \frac{n_k - 1}{n_k} b_k
\end{equation}
Notice that $\alpha_k \in ]-1 + 2/n_k, 1[$. %\frac{2}{n_k}, 1[$. 
Then:
$$ \avemap(A) = \begin{pmatrix}
   \alpha_1  & c_{1,2} & \dots & c_{1,p} \\
    c_{1,2} & \alpha_2 & \ddots & \vdots \\
    \vdots & \ddots  &  \ddots & c_{p-1,p}\\
    c_{1, p} & \dots & c_{p-1,p} & \alpha_{p} \\
    \end{pmatrix}. $$%

\subsection{The block filling map}
We call \textit{block filling map} from  $\mat{p}$ to $\blockClass{p}$, 
the linear map  
$$\fillmap(\m{C}) := 
  \begin{pmatrix} 
   C_{1,1} \, \m{J}_{n_1, n_1}  & \dots & C_{1,p} \,\m{J}_{n_1, n_p} \\
   \vdots & & \vdots  \\         
   C_{p,1} \, \m{J}_{n_p, n_1} & \dots & C_{p,p} \,\m{J}_{n_p, n_p} \\      
\end{pmatrix}.
$$
This transformation fills each block $(i,j)$ by a constant value, given by $C_{i,j}$.   
Notice that $\fillmap$ is a positive linear map, since
$ \fillmap(C) = \m{W}^\top \m{C} \m{W} $
where $W$ is the $p \times n$ matrix
$$ \m{W} := 
\begin{pmatrix} 
  \one{n_1}^\top & 0 & \dots & 0 \\
  0 & \one{n_2}^\top & \ddots & \vdots \\
  \vdots & \ddots & \ddots & 0\\
  0 & \dots & 0 & \one{n_p}^\top\\
\end{pmatrix}. 
$$

Clearly the filling and averaging operations cancel each other when they are done in this order, 
corresponding to $WV = I_p$ and to the relation
$$ \avemap(\fillmap(C)) = C.$$
On the other hand, $\fillmap(\avemap(C))$ is generally not equal to $C$, 
and $VW$ is an orthogonal projection of rank~$p$.

\subsection{A wider class of block matrices}
The results of the next section are valid on a wider class of block matrices, that we now introduce.
This is the set $\blockClassG{p}$ of symmetric block matrices  %  \supset \blockClass{p}$ 
$$ A = \begin{pmatrix}
    A_{1,1}  & \dots & A_{1,p} \\
    \vdots & & \vdots \\
	A_{p,1}  & \dots & A_{p,p}  \\
    \end{pmatrix}$$
where off-diagonal blocks $A_{k, \ell}$ are constant matrices of size $n_k \times n_\ell$,
and diagonal blocks $A_{k,k}$ are of size $n_k$ and verify $\m{A}_{k,k} - \mean{\m{A}_{k,k}} J_{n_k} \psdsymb 0$.\\

To see that $\blockClass{p} \subset \blockClassG{p}$, 
%let $A \in \blockClass{p}$ and denote $C = A_{k,k}$ for a given $k \in\intInt{p}$.
notice that %if $C$ is a p.d. CS matrix of size $n$ then $C - \mean{C} J_n \psdsymb 0$.  
if $C = \compound{n}{1, r}$ then 
$$\mean{C} = n^{-1}(1 + (n-1)r)$$ and
$$ C - \mean{C} J_n = (1-r) (I_n - n^{-1} J_n)$$ 
which is \pd{} for $r < 1$ because the symmetric matrix $\m{P} := \m{I_n} -n^{-1}\m{J_n}$ 
is an orthogonal projection: $\m{P} = \m{P}^2 = \m{P}\m{P}^\top$.\\

\subsection{Matrix pencil and positive definite interval}
When a parametric matrix of $\mat{p}$ depends linearly on a single parameter $t \in \R$, it can be written as a \emph{matrix pencil},
$$ t \mapsto D + tE $$
where $D, E \in \mat{p}$. 
The problem of finding $t$ such that $D + tE$ is positive definite was investigated in several studies (\cite{Caron_1986}, \cite{Caron_Traynor_Jibrin_2010}). 
A main result is that the admissible values form an open interval, 
called \emph{positive definite interval} (PD interval). 
Indeed it is the intersection of a straight line with the open convex cone $\mathcal{C}$ 
of positive definite matrices (\cite{Caron_1986}).
Now denote $d(c) = \det(D + c E)$.
When $D$ is \pd, $d$ is a continuous function such that $d(0) = \det(D) > 0$. % as the $\alpha_k$'s are $> 0$.
Furthermore, $d$ is null on the boundary of $\mathcal{C}$: 
if $A$ is \psd {} but not \pd, there exists $x \neq 0$ such that $x^\top A x = 0$
implying that $0$ is an eigenvalue of $A$.
The PD interval is then obtained by computing the negative and positive roots of $d$ which are closest to 0. 
Equivalently, the PD interval is related to the largest and smallest eigenvalues 
of the symmetric matrix $S = D^{-1/2} E D^{-1/2}$ (\cite{Caron_Traynor_Jibrin_2010}).
Indeed $d(c) = 0$ if and only if 
$$ \det \left(- \frac{1}{c} I_p - D^{-1} E \right) = 0$$
i.e. $-\frac{1}{c}$ is an eigenvalue of $D^{-1} E$, which has the same eigenvalues than $S$.\\

\section{Main results} \label{sec:main}
\subsection{Positive definiteness of elements of $\blockClass{p}$}
Positive definiteness of a block matrix (of size $n$) in $\blockClassG{p}$ 
-- and thus in $\blockClass{p}$ --   
can be simply checked on the block average matrix (of size $p$): 

\newtheorem{thm1}{Theorem}
\begin{thm1}\label{prop:TheoG}
Let $A \in \blockClassG{p}$.  
Then,
  \begin{enumerate}
  \item $\m{A} \psdsymb 0  \iff \avemap(A) \psdsymb 0$ and $\forall k \in \intInt{p}: \m{A}_{k,k} \psdsymb 0$.
  \item $\m{A} \pdsymb 0  \iff \avemap(A) \pdsymb 0$ and $\forall k \in \intInt{p}: \m{A}_{k,k} \pdsymb 0$.
  \end{enumerate}
\end{thm1}

The text of the theorem can be simplified on $\blockClass{p}$, for which diagonal blocks are CS correlation matrices.
Indeed, p.d. of $\m{A}_{k,k}$ is then equivalent to the two conditions $b_k < 1$ and $-(n_k - 1)^{-1} < b_k$
(Eq. \ref{eq:CompoudCarac}).
The former is contained in the definition of $\blockClass{p}$.
The latter is a consequence of p.d. of $\avemap(A)$ by looking at the signs of its diagonal elements (Eq. \ref{eq:alpha}). 
Hence, we have:

\newtheorem{thm0}[thm1]{Theorem}
\begin{thm0}\label{prop:averageMapping}
Let $A \in \blockClass{p}$. Then,%. Recall that $-1 < b_k < 1$ for all $k \in \intInt{p}$. Then,
$$A \pdsymb 0 \quad \iff \quad \avemap(A) \pdsymb 0 $$ 
\end{thm0}
\bigskip

For completeness, we provide a description of $\blockClassG{p}$, via a representation of its block diagonal elements.
%Such diagonal elements are connected with centered covariance matrices, which appear in the study of Gaussian process regression (see e.g. \cite{GRD_degeneracy}, Example 5).
It is shown that $\blockClassG{p}$ is in general much wider than $\blockClass{p}$:
the $k$-th diagonal block of $\blockClassG{p}$ is represented by $1 + n_k(n_k-1)/2$ parameters, 
whereas the $k$-th diagonal block of $\blockClass{p}$ is a CS correlation matrix described by $1$ parameter.\\
\newtheorem{prop1}{Proposition}
\begin{prop1}\label{prop:WiderClassCarac}
Let $R$ be a $n\times(n-1)$ matrix whose columns form a basis of $\one{n}^\perp$.
Then all $C$ symmetric matrices of size $n$ such that $\m{C} - \mean{\m{C}} J_{n} \psdsymb 0$
are written in a unique way as 
\begin{equation*}
  \label{eq:Repr}
  \m{C} = \mu J_n + \m{R} \m{B} \m{R}^\top
\end{equation*}
where $B$ is a \psd matrix of size $n-1$ and $\mu$ is a real number.
\end{prop1}
\bigskip

\subsection{Positive definiteness of elements of $\blockClassI{p}$}
We now focus on the subclass $\blockClassI{p}$, for which the positive definiteness condition 
of Theorem~\ref{prop:averageMapping} can be further simplified.
Let $A \in \blockClassI{p}$ and denote $B = \phi(A)$:
$$ B = \begin{pmatrix}
   \alpha_1  & c & \dots & c \\
    c & \alpha_2 & \ddots & \vdots \\
    \vdots & \ddots  &  \ddots & c \\
    c & \dots & c & \alpha_p\\
    \end{pmatrix} $$
where the $\alpha_k$'s are given by Eq.~\ref{eq:alpha} and lie in $]-1 + 2/n_k, 1[$. %\frac{2}{n_k}, 1[$.
We can further restrict the search of a necessary and sufficient condition to vectors 
$\bs{\alpha}$ such that $\alpha_k > 0$ ($k = 1, \dots, p$).
Indeed this condition is necessary as all diagonal minors must be strictly positive.
Notice that ``$0 < \alpha_k < 1$'' is the condition for positive definiteness of 
the block diagonals terms $\compound{n_k}{1,b_k}$ in $A$.\\

When $\bs{\alpha}$ is fixed, the map $c \mapsto B$ is a matrix pencil 
$$ B = D + c E$$
with $D = \diag(\alpha_1, \dots, \alpha_p)$ and $E = J_p - I_p$.
As recalled in Section~\ref{sec:preliminaries}, the values of $c$ such that $B$ is positive definite form an open interval, 
called positive definite interval (PD interval). 
It is obtained by computing the negative and positive roots of $d: c \mapsto \det(B)$ which are closest to 0.

In this context, our main contribution is to provide a localization of the roots of $d$ 
as well as an analytical expression. They both simplify the computation of the PD interval.

\newtheorem{lemma1}{Lemma}
\begin{lemma1}\label{prop:detLemma}
Let $B = \begin{pmatrix}
   \alpha_1  & c & \dots & c \\
    c & \alpha_2 & \ddots & \vdots \\
    \vdots & \ddots  &  \ddots & c \\
    c & \dots & c & \alpha_p\\
    \end{pmatrix}$, 
and denote by $\alpha_{(1)} \leq \dots \leq \alpha_{(p)}$ the values of $\alpha_1, \dots, \alpha_p$ rearranged in ascending order.\\
Denote $d_p(c) = \det(B)$ and more generally let $d_k(c)$ be the leading minor of $B$ ($k = 1, \dots, p$). 
Then $d_k$ is a polynomial of order $k$ given by
\begin{equation} \label{eq:detB}
d_k(c) = \prod_{m=1}^k (\alpha_m - c) + c \sum_{m = 1}^k \left( \prod_{\substack{1 \leq \ell \leq k \\ \ell \neq m}} (\alpha_l - c) \right).
\end{equation}
Let us further assume that all $\alpha_k's$ are $>0$. 
Then the roots $r_{k,1}, \dots, r_{k,k}$ of $d_k$ are all real and interlaced with the $\alpha_{(k)}$'s:
$$ r_{k,1} < 0 < \alpha_{(1)} \leq r_{k,2} \leq \alpha_{(2)} \leq \dots \leq r_{k,k} \leq \alpha_{(k)} < 1$$
Furthermore,
$$- \sqrt{\alpha_{(1)} \alpha_{(2)}} = r_{2,1} \leq \dots \leq r_{k, 1} < 0 < r_{k,2} \leq \dots \leq r_{2, 2} = \sqrt{\alpha_{(1)} \alpha_{(2)}}$$ 
\end{lemma1}
\bigskip
\newtheorem{thm2}[thm1]{Theorem}
\begin{thm2}\label{prop:CNSiso}
Let $A \in \blockClassI{p}$. With the notations of Lemma~\ref{prop:detLemma},
$$A \pdsymb 0 \quad \iff \quad 0 < \alpha_k < 1 \; \, (k \in \intInt{p}) \quad \textrm{and} \quad r_{p,1} < c < r_{p,2}$$
Furthermore, for all $k \in \intInt{p}$, the interval $]r_{p,1}, r_{p,2}[$ is increasing with $\alpha_k$.
\end{thm2}
\bigskip
Notice that the interval on $c$ given in Theorem~\ref{prop:CNSiso} can be easily found numerically. 
Indeed, from Lemma~\ref{prop:detLemma}, $r_{p,1}$  and $r_{p,2}$  are roots of $c \mapsto \det(\phi(A))$ that are perfectly localized,
and can be found by a zero search algorithm such as Brent's algorithm (\cite{Brent_1973_book}). 
More precisely $r_{p,1}$ is the unique root in the interval $]-\sqrt{\alpha_{(1)}\alpha_{(2)} }, 0[$,
and $r_{p,2}$ the unique root in $]0, \sqrt{\alpha_{(1)}\alpha_{(2)} }[$.

\section{Corollaries} \label{sec:cor}
We first give a direct consequence of Theorem~\ref{prop:averageMapping} for two groups, 
and a sufficient condition for elements of $\blockClassI{p}$ obtained by comparison to a compound symmetry matrix. 
They are expressed with the $\alpha_k$'s of Eq.~\ref{eq:alpha}.

\newtheorem{cor1}{Corollary}
\begin{cor1}[Case of 2 groups]\label{prop:2groups}
Let $A \in \blockClass{2}$. Then, 
$$A \pdsymb 0  \iff  0 < \alpha_k < 1 \;\, (k = 1,2) \quad \textrm{and} \quad \vert c_{1,2} \vert < \sqrt{\alpha_1 \alpha_2} $$
\end{cor1}

\newtheorem{cor2}[cor1]{Corollary}
\begin{cor2}[A sufficient condition on $\blockClassI{p}$]\label{prop:CS}
Let $A \in \blockClassI{p}$ and $\alpha^\star = \min_{1 \leq k \leq p}(\alpha_k)$. 
If $0 < \alpha_k < 1$ ($k \in \llbracket 1, p \rrbracket$) and $- \frac{\alpha^\star}{p-1} < c < \alpha^\star$,
then $A \pdsymb 0$.  
\end{cor2}

For completeness, we give a specific result when all groups have the same size, 
a common between group correlation, and a common within group correlation.
There is nothing new here, as it is a special case of block CS matrices for which the result is known,
but this gives another view of it.

\newtheorem{cor3}[cor1]{Corollary}
\begin{cor3}\label{prop:toeplitzLike}
Assume that $n_1 = \dots = n_p = n_0 := n/p$ and for all $k \in \llbracket 1, p \rrbracket$, $b_k = b$. % \in ]-1, 1[$. 
Define $\alpha = (1 + (n_0-1)b)/n_0$. Then,  
$$A \pdsymb 0 \quad \iff \quad 0 < \alpha < 1 \quad \textrm{and} \quad - \frac{\alpha}{p-1} < c < \alpha.$$
\end{cor3}

\section{A numerical application} \label{sec:num}
\begin{figure}[!ht]
\centering
\includegraphics[width=0.95\linewidth]{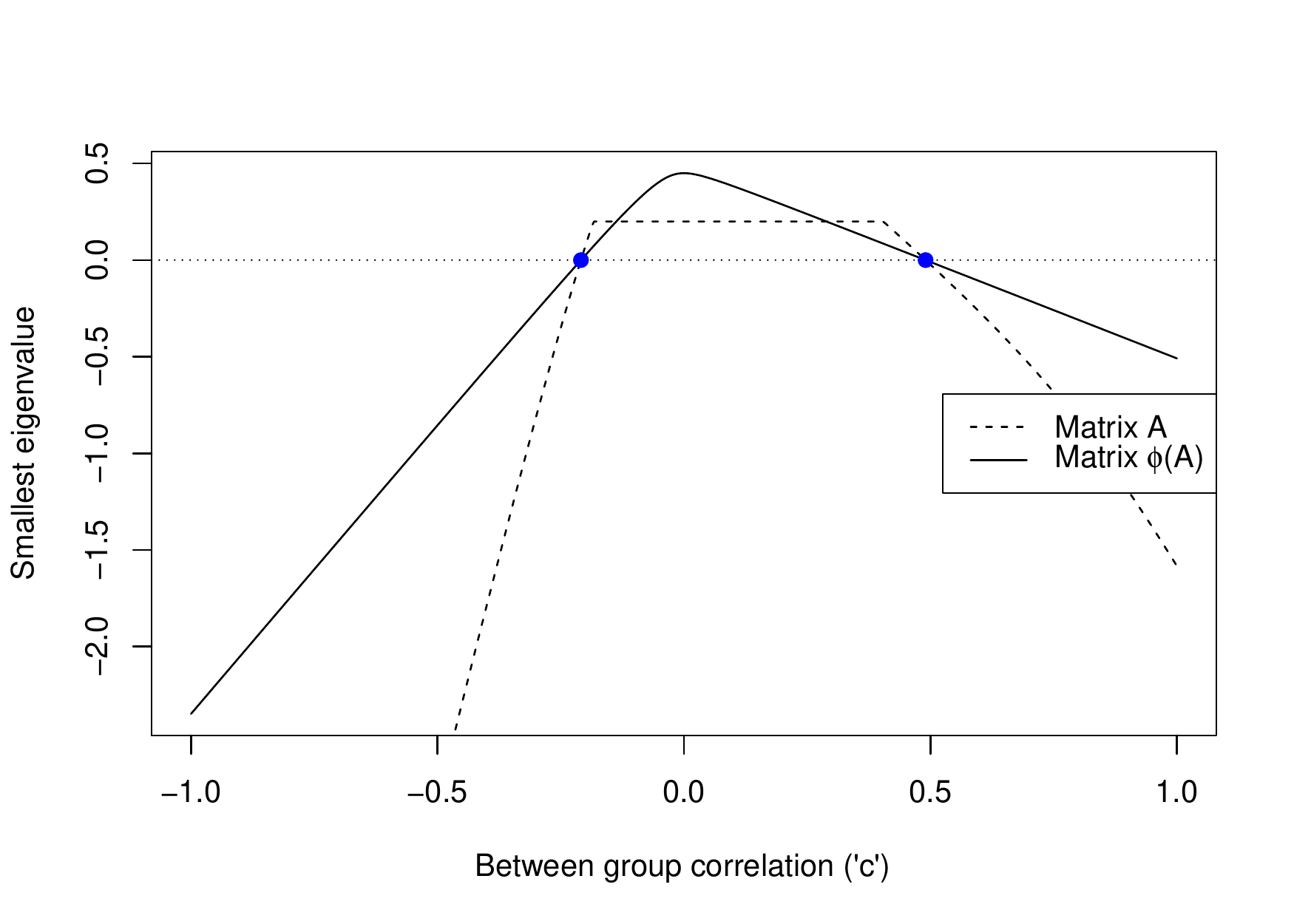} 
\caption{Smallest eigenvalue of the $15 \times 15$ block matrix A of the example
and the $4 \times 4$ block average matrix $\avemap(A)$ as a function of the between group correlation $c$.}\label{fig:app1}
\end{figure}

Consider the parametric matrix $A \in \blockClassI{15}$ corresponding to 4 groups of size $2, 4, 3, 6$ with a common between group correlation.
The condition $0 < \alpha_k < 1$ means that: 
$$ -1 < b_1 < 1, \quad -\frac{1}{3} < b_2 < 1, \quad -\frac{1}{2} < b_3 < 1, \quad  -\frac{1}{5} < b_4 < 1$$
Now choose $b_1 = -0.1, b_2 = 0.4, b_3 = 0.7, b_4 = 0.8$. 
The sufficient condition of Corollary~\ref{prop:CS} gives the following interval for the between group parameter
(printing is limited to the first 2 digits) :
$$ -0.03 < c < 0.45 $$
Theorem~\ref{prop:CNSiso} gives the optimal interval, which is a bit larger:
\begin{equation} \label{eq:examplePDInterval}
 -0.20 < c < 0.49 
\end{equation}
This interval was obtained by using the Brent's algorithm, implemented in R (function \texttt{uniroot} of package \texttt{stats}, 
\cite{R}), thanks to the localization given in Theorem~\ref{prop:CNSiso}.
The result is illustrated on Figure \ref{fig:app1}. 
For a pedagogical purpose, we represented the smallest eigenvalue of $A$ and $\avemap(A)$ as a function of $c$, 
denoted respectively $\lambda_{A}(c), \lambda_{\phi(A)}(c)$. These eigenvalues have been computed numerically using the \texttt{eigen} function in R.
This confirms the results of Theorem~\ref{prop:averageMapping} and Theorem~\ref{prop:CNSiso}:
the values of $c$ such that the smallest eigenvalue is positive is the same in both cases, 
and correspond to the theoretical interval (\ref{eq:examplePDInterval}).
Notice that the size of $\avemap(A)$ is much smaller than $A$. %, and may be preferred. 
Finally, the plateau on the graph $c \mapsto \lambda_{A}(c)$ is explained by the fact that $\lambda_{A}(0)$ is equal to
$1 - \max_{1 \leq k \leq p} {b_k}$, which does not depend on $c$, 
which remains true on an open interval containing $0$ by continuity of the determinant.

\section{Proofs} \label{sec:proofs}

\subsection{Proof of Theorem~\ref{prop:TheoG} and related material.}

\begin{proof}[Proof of Proposition~\ref{prop:WiderClassCarac}]
First of all, notice that finding a symmetric matrix $C$ of size $n$ such that $C - \mean{C} J_n \psdsymb 0$
is equivalent to finding a \psd~matrix $C_0$ of size $n$ such that $\mean{C_0} = 0$.
The link between $C$ and $C_0$ is given by $C = \mu J_n + C_0$, where $\mu$ is a real number equal to $\mean{C}$.\\
Hence, the proposition will be proved if we show that all \psd~matrices $C_0$ such that $\mean{C_0} = 0$
are of the form 
\begin{equation} \label{eq:RepresC}
\m{C_0}(\m{B}) = \m{R} \m{B} \m{R}^\top
\end{equation}
where $B$ is a \psd~matrix of size $n$,
and $\m{R}$ is a $n \times (n-1)$ matrix whose columns form a basis of $\one{n}^\perp$.\\ 
Notice that all \psd~matrix $\m{C_0}$ of the form~\ref{eq:RepresC} satisfy, by definition of $R$,
$$\mean{C_0(B)} = n^{-2} \one{n}^\top R \m{B} R^\top \one{n} = 0.$$
Conversely, let $C_0$ be a \psd~matrix such that $\mean{C_0}=0$. 
Thus $C_0$ is the covariance matrix of a random vector $Z$.
Without loss of generality, we can assume that $Z$ is centered: $E(Z)=0$.
Then, $\mean{Z} =  n^{-1} \one{n}^\top Z$ is also centered and 
$$\var(\mean{Z}) = n^{-2} \one{n}^\top C_0 \one{n} = \mean{C_0} = 0.$$
Hence, $\mean{Z} = 0$ with probability 1.
Thus $\one{n}^\top Z = 0$ and $Z$ is spanned by the columns of $R$:
there exists a random variable $U$ such that $Z  = R U$.
As a consequence, $\cov(Z) = R B R^\top$, with $B = \cov(U)$.\\

To see that $B$ and $\mu$ are uniquely determined, assume that $C = \mu' J_n + R B' R^\top$ 
holds for another \psd~matrix $B'$ of size $n-1$ and another real $\mu'$. 
Then $\mu' = \mu = \mean{C}$. Thus, $R B' R^\top = R B R^\top$.
By multiplying this relation at left side by $R^\top$ and by $R$ at right side, we get $B = B'$ 
because the Gram matrix $R^\top R$ has rank $n-1$ and is thus invertible. % = I_{n-1}$.
\end{proof}

We then give a technical lemma.

\newtheorem{lemma2}[lemma1]{Lemma}
\begin{lemma2} \label{lemma:l0}
Let $\m{C}$ be a $n \times n$ symmetric matrix 
%and $\bar{\m{C}} := \bar{C} \, \m{J}$ where $\bar{C}$ is the average element of $\m{C}$. 
such that $\m{C} - \mean{\m{C}}J_n$ is positive.  
If $\m{C}$ is p.d. and  if $\bs{\beta}^\top
    \left[\m{C} - \mean{\m{C}} J_n \right] \bs{\beta} = 0$ for some
    vector $\bs{\beta}$, then $\bs{\beta}$ has constant elements, i.e. $\bs{\beta} \propto \one{n}$.
\end{lemma2}

\begin{proof}[Proof of Lemma~\ref{lemma:l0}]
Since $\m{C}$ has rang $n$ and $\mean{\m{C}} J_n$ has rank~$1$, 
the matrix $\m{C_0} := \m{C} - \mean{\m{C}} J_n$ must have rank $n-1$. 
Thus, by Prop.~\ref{prop:WiderClassCarac}, $C_0 = RBR^\top$ where $B$ is p.d. of size $n-1$,
and $R$ is a $n \times (n-1)$ matrix whose columns form an orthonormal basis of $\one{n}^\perp$.
Now, $0 = \bs{\beta}^\top \m{C_0} \bs{\beta} = \bs{\beta}^\top R B R^\top \bs{\beta}$.
Since $B$ is p.d. it implies $\bs{\beta}^\top R = 0$. 
By definition of $R$, we have $\bs{\beta} \propto \one{n}$.
\end{proof}

\begin{proof}[Proof of Theorem~\ref{prop:TheoG}]
\medskip\par\noindent  \textit{Statement 1.} 
If $\m{A}$ is positive then the diagonal blocks $\m{A}_{k,k}$ are also positive, 
and $\avemap(A)$ is positive since $\avemap$ is a positive linear map.\\
Conversely, if $\avemap(A)$ is positive then $\fillmap(\avemap(A))$ is positive, 
as $\fillmap$ is a positive linear map. 
Now, since $A$ has constant off-diagonal blocks, 
$A - \fillmap(\avemap(A))$ is a block diagonal matrix with diagonal blocks 
$\m{A}_{k,k} - \mean{\m{A}_{k,k}}$, which are assumed to be positive. 
Hence $\m{A}$ is positive as a sum of two positive matrices. 
 
\medskip\par\noindent\textit{Statement 2.} If $\m{A}$ is p.d. then
  so are the $\m{A}_{k,k}$. Moreover $\fillmap(A) = W^\top A W$ is p.d.
  since the $p \times n$ matrix $\m{W}$ has rank $p$.\\  
  Now assume that $\avemap(\m{A})$ and all the $\m{A}_{k,k}$'s are p.d. 
  Denote $\m{A}^\star := \fillmap(\avemap(\m{A}))$, 
  and let $\bs{\beta}$ be a vector of length $n$ with $p$ component vectors $\bs{\beta}_k$. 
  Then,
  $$
  \bs{\beta}^\top \m{A} \bs{\beta}  =  
  \bs{\beta}^\top \left\{ \m{A}  -\m{A}^\star \rule{0pt}{0.9em}
  \right\} \bs{\beta}
  +  \bs{\beta}^\top \m{A}^\star  \bs{\beta} =: u + v.
  $$ 
  By Statement 1, both matrices $\m{A} -\m{A}^\star$ and $\m{A}^\star$ are
  positive. Hence, the left hand side is equal to zero iff both terms at the
  right hand side vanish. Now $\m{A} -\m{A}^\star$ is %the matrix between the curly brackets is
  block-diagonal, and each diagonal block is the difference between a
  block $\m{A}_{k,k}$ and its average element. By lemma~\ref{lemma:l0}, we
  have $u=0$ iff each vector $\bs{\beta}_k$ is constant ($k \in \intInt{p}$). % to  $p$. 
  The second term $v$ is a quadratic form in the vector of group sums 
  $\bs{\gamma} := W \bs{\beta} = [n_1 \mean{\beta_1}, \dots, n_p \mean{\beta_p}]^\top$, 
  namely $v = \bs{\gamma}^\top \avemap(\m{A}) \bs{\gamma}$. 
  So $v$ vanishes iff each average $\mean{\beta_k}$ does. 
  Thus if $\bs{\beta}^\top \m{A} \bs{\beta} =0$, each block $\bs{\beta}_k$ is constant 
  and has zero mean, which is only possible when $\bs{\beta} = \m{0}$.
\end{proof}

\subsection{Proof of Theorem~\ref{prop:CNSiso}} %Lemma~\ref{prop:detLemma}}
\begin{proof}[Proof of Lemma~\ref{prop:detLemma}]
For the sake of simplicity, and without restriction, 
we can assume that the $\alpha_k's$ have been sorted in increasing order: 
$\alpha_{(k)} = \alpha_k$ ($k = 1, \dots, p)$.
Indeed the reordering does not change the determinant.\\
Now for $k = 1, \dots, p$, define: 
$$d_k(c) = \det (B_{[1:k, 1:k]})$$
where $B_{[1:k, 1:k]}$ is the matrix extracted from $B$ by keeping the first $k$ rows and first $k$ columns.
In particular $d_p(c) = d(c) = \det(B)$, and $d_2(c) = \alpha_1 \alpha_2 - c^2$.\\

Let us now derive the analytical expression of $\det(B)$. 
We first prove a recurrence formula.
By expanding with respect to the last column, we have:
$$ d_k(c) = \det \begin{pmatrix}
   \alpha_1  & c & \dots &  c & 0\\
    c & \ddots & \ddots &  \vdots & \vdots \\
    \vdots & \ddots & \ddots  &  c & \vdots\\
    c & \dots & c & \alpha_{k-1}  & 0\\
    c & \dots & \dots & c & \alpha_{k} - c \\
    \end{pmatrix}
+
   \det \begin{pmatrix}
   \alpha_1  & c & \dots &  c & c\\
    c & \ddots & \ddots &  \vdots & \vdots \\
    \vdots & \ddots & \ddots  &  c & \vdots\\
    c & \dots & c & \alpha_{k-1}  & c\\
    c & \dots & \dots & c & c \\
 \end{pmatrix}    $$ 
The first term is equal to $(\alpha_k - c)d_{k-1}(c)$. 
For the second one, subtracting the last column to the other ones gives:
$$ \det \begin{pmatrix}
   \alpha_1  & c & \dots &  c & c\\
    c & \ddots & \ddots &  \vdots & \vdots \\
    \vdots & \ddots & \ddots  &  c & \vdots\\
    c & \dots & c & \alpha_{k-1}  & c\\
    c & \dots & \dots & c & c \\
 \end{pmatrix}     =
 \det \begin{pmatrix}
   \alpha_1 - c  & 0 & \dots & 0 & c\\
    0 & \ddots & \ddots &  \vdots & \vdots \\
    \vdots & \ddots & \ddots  &  0 & \vdots\\
    0 & \dots & 0 & \alpha_{k-1} - c  & c\\
    0 & \dots & \dots & 0 & c \\
 \end{pmatrix}     $$
Finally we obtain a relation linking $d_k$ to $d_{k-1}$, valid for $k \geq 2$:
\begin{equation} \label{eq:detBlocRecurrence}
d_k(c) = (\alpha_k - c) d_{k-1}(c) + c  \prod_{m = 1}^{k-1} (\alpha_m - c)
\end{equation} 
The explicit formula of $d_k(c)$ (written in Eq.~\ref{eq:detB} for $k=p$) is then proved recursively on $k$. 
Indeed, one can easily check it for $k = 2$, 
and if it is true for $k-1 \geq 2$, then with Eq.~\ref{eq:detBlocRecurrence}, it holds:
\begin{eqnarray*}
d_k(c) &=& (\alpha_k - c) \prod_{m=1}^{k-1}(\alpha_m - c) 
+ c \sum_{m = 1}^{k-1} \left( \prod_{\substack{1 \leq \ell \leq k-1 \\ \ell \neq m}} (\alpha_\ell - c) (\alpha_k - c) 
+ \prod_{m = 1}^{k-1} (\alpha_m - c) \right )  \\
&=& \prod_{m=1}^k (\alpha_m - c) + c \sum_{m = 1}^k \left( \prod_{\substack{1 \leq \ell \leq k \\ \ell \neq m}} (\alpha_\ell - c) \right) 
\end{eqnarray*}
We also deduce of Eq.\ref{eq:detBlocRecurrence} that $d_k$ is a polynomial of order $k$ whose leading term is $(-1)^{k-1} (k-1) c^k$. In particular, $d_k(c) \rightarrow - \infty $ when $c \rightarrow - \infty$.\\% is $-\infty$.\\

From now on, assume that $0 < \alpha_1 < \dots < \alpha_p$.
The interlacing of the roots of $d_k$ with the $\alpha_k$'s is obtained by evaluating $d_k$ at $\alpha_\ell$ ($1 \leq \ell \leq k$).
Indeed, subtracting column $\ell$, whose components are all equal to $\alpha_\ell$, one directly obtains:
$$ d_k(\alpha_\ell) = \alpha_\ell \prod_{m \neq \ell} (\alpha_m - \alpha_\ell) $$
Thus the sign of $d_k(\alpha_\ell)$ depends on the rank of $\ell$ in $\{1, \dots, k\}$. 
Those signs alternate when one read the sequence $\{1, \dots, k\}$ starting by $+1$ at $\ell = 1$.
Since the limit in $- \infty$ is $-\infty$, and as the number of zeros of $d_k$ is less than $k$, 
we deduce by the intermediate value theorem that there is exactly 1 zero in each interval 
$]-\infty, \alpha_1], [\alpha_1, \alpha_2], \dots, [\alpha_{k-1}, \alpha_k]$.
Furthermore $d_k(0) = \alpha_1 \dots \alpha_k > 0$, resulting in the interlacing 
$$ r_{k,1} < 0 < \alpha_1 < r_{k,2} < \alpha_2 < \dots < r_{k,k} < \alpha_k < 1$$

We are now able to locate the roots of $d_k$ relatively to those of $d_{k-1}$.
From Eq~\ref{eq:detBlocRecurrence}, it holds:
$$ d_k(r_{k-1,\ell}) = r_{k-1,\ell} \prod_{m = 1}^{k-1}  (\alpha_m - r_{k-1,\ell})$$
In particular $d_k(r_{k-1,1}) < 0$, showing that $r_{k-1,1} < r_{k, 1} < 0$
since $d_k$ is negative on $]-\infty, r_{k,1}]$ and positive on $[r_{k,1}, \alpha_1]$.
Similarly, one deduce from Eq.~\ref{eq:detBlocRecurrence} that $d_k(r_{k-1,2}) < 0$ 
so that $ 0 < r_{k,2} < r_{k-1, 2}$ since $d_k$ is positive on $[\alpha_1, r_{k,2}]$ and negative on $[r_{k,2}, \alpha_2]$.
Finally, we obtain by recursion that
$$- \sqrt{\alpha_1 \alpha_2} = r_{2,1} < \dots < r_{k, 1} < 0 < r_{k,2} < \dots < r_{2, 2} = \sqrt{\alpha_1 \alpha_2}$$ 
More generally, Eq.~\ref{eq:detBlocRecurrence} implies that 
the positive roots of $d_k$ are interlaced with the positive roots of $d_{k-1}$. 
Finally one can check that the previous results remain valid in case of equality of several $\alpha_k's$,
where the corresponding strict inequalities are replaced by large inequalities.\\

For completeness, we give an alternative proof of this last statement, 
using algebraic arguments. 
Consider the pencil matrix $c \mapsto B = D + c E$.
Recall that $c$ is a root of $d_p$ if and only if 
%Following \cite{Caron_Traynor_Jibrin_2010}, notice that $d_p(c) = 0$ if and only if 
$$ \det \left(- \frac{1}{c} I_p - D^{-1} E \right) = 0$$
i.e. $-\frac{1}{c}$ is an eigenvalue of $D^{-1} E$, 
or equivalently of the symmetric matrix $S = D^{-1/2} E D^{-1/2}$.
Similarly, $c$ is the root of $d_k$ if and only if 
$-\frac{1}{c}$ is an eigenvalue of 
$$D_{[1:k, 1:k]}^{-1/2} E_{[1:k, 1:k]} D_{[1:k, 1:k]}^{-1/2} = S_{[1:k, 1:k]},$$
where we use the fact that $D$ is diagonal for the last equality. 
Hence the roots of $d_k$ are the negative inverse of the eigenvalues 
of the $k^{\mathrm{th}}$ principal submatrix of $S$.
Furthermore, by Sylvester's law, $S_{[1:k, 1:k]}$ has the same inertia as $E_{[1:k, 1:k]} = J_k - I_k$.
It is easy to see that $E_{[1:k, 1:k]}$ has two eigenvalues: $k-1$ with multiplicity $1$ (eigenvector $\one{k}$) 
and $-1$ with multiplicity $k-1$ (eigenspace: $\one{k}^\perp$). 
Thus $S_{[1:k, 1:k]}$ has $1$ positive eigenvalue, and $k-1$ negative ones. 
Consequently, $d_k$ has $1$ negative root and $k-1$ positive ones.
Now by Cauchy interlacing theorem (\cite{Horn_Johnson_book}, \S 4.3.), 
the eigenvalues of consecutive principal submatrices are interlaced.
As $c \mapsto -\frac{1}{c}$ is increasing on $\mathbb{R}_-^*$ and $\mathbb{R}_+^{*}$,
it results in two separate interlacings for the roots of consecutive $d_k$'s: one for negative roots, one for positive roots. 
We finally get the announced inequalities:
$$- \sqrt{\alpha_{(1)} \alpha_{(2)}} = r_{2,1} \leq \dots \leq r_{k, 1} < 0 < r_{k,2} \leq \dots \leq r_{2, 2} = \sqrt{\alpha_{(1)} \alpha_{(2)}}$$ 
\end{proof}

\begin{proof}[Proof of Theorem~\ref{prop:CNSiso}]
The result is straightforward from Lemma~\ref{prop:detLemma}, 
implying that the PD interval on $c$ is given by $]r_{p,1}, r_{p,2}[$.
An alternative proof is given by applying the Sylvester conditions and using Lemma~\ref{prop:detLemma}.
Indeed %when $\bm{\alpha}$ is fixed, 
the intersection of sets $\{c, d_k(c) > 0 \}$ when $k \in \llbracket 1, p \rrbracket$  
is the interval $]r_{p,1}, r_{p,2}[$.\\
Furthermore, let $k \in \llbracket 1, p \rrbracket$, and $0 < \alpha_k < \alpha_k'$.
Denote by $D'$ the diagonal matrix obtained by replacing $\alpha_k$ by $\alpha_k'$ in $D$.
Then $D' + cE > D + cE$, which shows that the PD interval is increasing with $\alpha_k$.
\end{proof}

\subsection{Proof of corollaries}
\begin{proof}[Proof of Corollary~\ref{prop:CS}]
Let $B^\star = \compound{p}{\alpha^\star, c}$. 
By definition of $\alpha^\star$, $\phi(A) \psdsymb B^\star$. 
Now $B^\star \pdsymb 0 \iff - \frac{\alpha^\star}{p-1} < c < \alpha^\star$.
The result follows from Theorem~\ref{prop:averageMapping}.
\end{proof}

\begin{proof}[Proof of Corollary~\ref{prop:toeplitzLike}]
Here $\avemap(A) = \compound{p}{\alpha, c}$ which gives the announced necessary and sufficient condition 
by Theorem~\ref{prop:averageMapping}.\\ 
Alternatively the result can be obtained directly on $A$, 
which is here a $p \times p$ block compound symmetry covariance matrix
$$ A = \begin{pmatrix}
    \Sigma_0  & \Sigma_1 & \dots &  \Sigma_1 \\
    \Sigma_1 & \ddots & \ddots & \vdots \\
    \vdots & \ddots & \ddots &  \Sigma_1\\
    \Sigma_1 & \dots & \Sigma_1 & \Sigma_0\\
    \end{pmatrix} $$
with $\Sigma_0 = \compound{n_0}{1,b}$ and $\Sigma_1 = c J_{n_0}$.
It is known (see e.g. \cite{Ritter_Gallegos_2002}, Lemma 4.3.) that $A \pdsymb 0$ if and only if $- \frac{\Sigma_0}{p-1} < \Sigma_1 <  \Sigma_0$.
Now, $\Sigma_0 - \Sigma_1 = \compound{n_0}{1-c, b-c}$ 
is \pd{} if and only if $ - \frac{1-c}{n_0-1} < b-c < 1-c$ leading to $c < \alpha$.
Similarly, $\Sigma_1 +  \frac{\Sigma_0}{p-1} = \compound{n_0}{c+\frac{1}{p-1}, c+\frac{b}{p-1}}$ 
is \pd{} if and only if $ - \frac{c+\frac{1}{p-1}}{n_0-1} < c+\frac{b}{p-1} < c+\frac{1}{p-1}$ 
leading to  $- \frac{\alpha}{p-1} < c$.
\end{proof}

\section*{Acknowledgements}
Part of this research was conducted within the frame of the Chair in Applied Mathematics OQUAIDO, gathering partners in technological research (BRGM, CEA, IFPEN, IRSN, Safran, Storengy) and academia (CNRS, Ecole Centrale de Lyon, Mines Saint-Etienne, University of Grenoble, University of Nice, University of Toulouse) around advanced methods for Computer Experiments.

\nocite*

%\bibliographystyle{plain}
%\bibliography{blockMatrices}

\end{document}